\documentclass[letterpaper,reqno, 10 pt, conference]{ieeeconf} 
\IEEEoverridecommandlockouts
\usepackage[utf8]{inputenc}
\usepackage[english]{babel}
\usepackage{float}
\usepackage{amsmath}
\usepackage{amsfonts}
\usepackage{amssymb}
\usepackage{verbatim}
\usepackage{algpseudocode}
\usepackage{algorithm}
\usepackage{array}
\usepackage{subcaption}
\usepackage{cite}
\usepackage{tikz}
\usepackage{tabularx,multirow,booktabs,blindtext}
\usepackage{graphicx}
\graphicspath{{images/}}
\usepackage[export]{adjustbox}

\usepackage{amsthm}
\theoremstyle{plain}
\newtheorem{thm}{Theorem}
\newtheorem{lem}{Lemma}[section]
\newtheorem{prop}{Proposition}[section]

\theoremstyle{definition}

\newtheorem{prbm}{Problem}

\newtheorem{assum}{Assumption}
\theoremstyle{remark}
\newtheorem{rem}{Remark}

\usetikzlibrary{automata,arrows,calc,positioning}

\def\R{\mathbb{R}}
\def\C{\mathbb{C}}

\usepackage{footnotehyper} \makesavenoteenv{algorithm}

\begin{document}

\title{Privacy-preserving cloud computation of algebraic Riccati equations }

\author{M.V. Surya Prakash, Nima Monshizadeh
\thanks{The authors are with the Engineering and Technology Institute, University of Groningen, The Netherlands. Emails: {\tt\small v.s.p.malladi@rug.nl}, {\tt\small n.monshizadeh@rug.nl}.
The work of authors is supported by NWO grant OCENW.KLEIN.257.}}

\maketitle

\begin{abstract}
We address the problem of securely outsourcing the solution of algebraic Riccati equations (ARE) to a cloud. Our proposed method explores a middle ground between privacy preserving algebraic transformations and perturbation techniques, aiming to achieve simplicity of the former and strong guarantees of the latter. Specifically, we modify the coefficients of the ARE  
in such a way that the cloud computation on the modified ARE returns the same solution as the original one, which can be then readily used for control purposes. Notably, the approach obviates the need for any algebraic decoding step. We present privacy-preserving algorithms with and without a realizability requirement, which asks for
preserving sign-definiteness of certain ARE coefficients in the modified ARE.  
For the LQR problem, this amounts to ensuring that the modified ARE coefficients
can be realized again as an LQR problem  for a (dummy) linear system. 
The algorithm and its computational load is illustrated through a numerical example. 
\end{abstract}

\section{Introduction}\label{sec:Introduction}

The rapid progress in science and technology in our highly interconnected society has led to an enormous volume of user data being gathered and analyzed through cyber networks. Nevertheless, these technological advancements have had a substantial impact on the privacy of individuals in society. Depending on their capabilities, curious parties can deduce sensitive and private information about users or systems from data from unsecured networks or servers.

The problem of preserving privacy of a cyber-physical system has been approached in multiple ways. The most popular ones are encryption based techniques, differential privacy, and information-theoretic methods. 

Various types of encryption schemes have been used to implement encrypted controllers in control over the cloud \cite{farokhi2016secure, kim2016encrypting} and distributed computation/control \cite{zhang2018enabling,ruan2019secure,lu2018privacy, hosseinalizadeh2022private} settings; see \cite{darup2021encrypted} for a brief introduction to the use of encrypted controllers in different classes of systems. A labeled homomorphic encryption scheme is employed to privately outsource LQG to a cloud in \cite{alexandru2019encrypted}. The encryption methods are highly secure, but this comes at the cost of very high computational load.   

Moreover, differential privacy based methods have been devised
for estimation \cite{le2013differentially}, distributed control \cite{huang2014cost} and distributed
convex optimization \cite{han2016differentially}, where data/models are perturbed by noise as a mean to hide privacy-sensitive information. Privacy-preserving Gaussian mechanisms have been proposed for stochastic dynamical systems in \cite{hayati2021finite,hayati2023infinite}.
While noise perturbation techniques are applicable to a wide range of problems, they generally suffer from performance loss due to the noise addition. 

Algebraic transformation approaches have been used for preserving privacy in optimization \cite{mangasarian2011privacy, vaidya2009secure}. These methods have found use in control theory as well, where, typically, sensitive data such as inputs and outputs and/or the system parameters are transformed before being communicated. For instance, in \cite{sultangazin2020symmetries}, linear maps are used to algebraically encode and decode system matrices and inputs-outputs for outsourcing an optimal control problem to the cloud. 
The work in \cite{monshizadeh2019plausible} proposes to augment the dynamics to hide privacy-sensitive parameters, e.g. using exchanged auxiliary variables that asymptotically track the state variables in the consensus protocol. Vanishing time-varying masks are proposed in \cite{altafini2020system} to obscure the private variables in multi-agents systems. 

In this work, we consider outsourcing the computation of solving an algebraic Riccati equation (ARE) to a cloud in a private fashion. Algebraic Riccati equations have found use in a wide range of disciplines such as control theory, fluid dynamics, game theory and transportation science \cite{chen1992necessary,bini2011numerical}. The ARE is a central equation in systems and control theory, with numerous applications including linear quadratic regulation (LQR), $H_2$ and $H_{\infty}$ control and (sub)optimal nonlinear control using State-Dependent Riccati Equation(SDRE)\cite{ccimen2008state}.  

The algorithms available to solve Riccati equations can be broadly divided into two categories: subspace algorithms and iterative algorithms. For a general system of size $n$, the time complexity of the fastest available algorithm is $\mathcal{O}(n^3)$ and the memory requirement is $\mathcal{O}(n^2)$ \cite{bai2000templates}. Therefore, finding solutions to the ARE becomes very expensive for large scale systems (e.g. $n > 1000$), which may exceed the local computational capacity. 
The issue is compounded when the dynamics are obtained via linearization from a nonlinear system whose operating point may change. The implementation of SDRE technique generally involves solving an ARE with updated coefficients at each sampling time \cite{dai2021nonlinear}. Furthermore, outsourcing the computation to the cloud frees up the local computational capacity which can be used for other purposes such as monitoring and fault detection. 

Outsourcing an ARE typically involves sending the ARE coefficients that typically correspond to sensitive system dynamics and optimization parameters; e.g. in the case of LQR, $H_2$ and $H_{\infty}$ optimal control. Communicating the true ARE coefficients to the cloud thus increases the risk of disclosure of plant parameters to malicious entities in a couple of ways: (i) The cloud can act as an adversary and reveal the acquired sensitive information to malicious third parties. (ii) Even if the cloud is trustworthy, it may fall victim to eavesdropping or cyber-attacks, thereby leaking the privacy-sensitive information. We note that several attacks on cyber-physical systems are enabled or aided by the knowledge of the system dynamics. For instance, the design and implementation of zero-dynamics attack and bias-injection attack are facilitated by the knowledge of system dynamics, whereas covert attacks can be launched by combining knowledge on the system dynamics with input-output measurements \cite{teixeira2015secure,pasqualetti2013attack}. We also emphasize that some of the system/optimization information can be inherently private; e.g. knowledge of the cost parameters in the corresponding optimal control problem.

The main idea put forward in this paper is to perturb the system model such that the result of the cloud's computation, namely the solution to the ARE, remains valid for the true system. The perturbation proposed here is inspired by the eigenvalue shifting technique in ARE theory \cite[Ch. 2]{bini2011numerical}.
Note that, unlike other algebraic techniques e.g. \cite{sultangazin2020symmetries,weeraddana2013per}, the proposed method does not require any decoding step to be applied to the solution computed by the cloud. As such, the proposed algebraic approach leverages the strong privacy offered by random/arbitrary perturbations without compromising the accuracy of the computation. 

The structure of the paper is as follows. Notations and preliminaries are provided in Section \ref{sec:Preliminaries}. The problems of interest are formulated in Section \ref{sec:Problem Statement}. Section \ref{sec:Main Results} presents the main results of the paper, namely providing algorithms that solve the formulated problems. The extent of the ambiguity faced by the cloud in retrieving the true system parameters is discussed in Section \ref{sec:Privacy analysis}. The application of the results to the notable case of linear quadratic regulators and a numerical illustration of the algorithms is provided in Section \ref{sec:Case Study}. The paper closes with concluding remarks in Section \ref{sec:Conclusions}. The proofs are collected in the appendix in Section \ref{sec:Appendix}.

\section{Notations and Preliminaries}\label{sec:Preliminaries}
The set of $n \times n$ real symmetric matrices is denoted by $\mathbb{S}_n$. The identity matrix and zero matrix of size $n \times n$ are denoted by $I_n$ and $\mathbf{0}_n$, respectively. The conjugate of a vector/matrix $A$ is denoted by $\bar{A}$ and its conjugate transpose by $A^*$. Given matrices $X$ and $Y$, we use the notation ${\rm col}(X,Y)$ to represent the concatenated matrix $\begin{bmatrix}
     X^T & Y^T
 \end{bmatrix}^T$.
 
The multispectrum of $A \in \C^{n \times n},$ denoted by $\sigma(A)$, is the multiset of eigenvalues of $A$. The size of a Jordan block with an eigenvalue $\lambda$ in the Jordan normal form of $A$ is called its partial multiplicity. We call $(\lambda,v)$ ($(\lambda,w)$) a right (left) eigenpair of $A$ if $\lambda$ is an eigenvalue of $A$ and $v$ ($w$) is the corresponding right (left) eigenvector.

A linear subspace $\mathcal{V}$ is called $A-$invariant if $A \mathcal{V} \subseteq \mathcal{V}$. If the columns of $V \in \R^{n \times m}$ form a basis for $\mathcal{V}$, then, there exists a nonsingular matrix $Z \in \R^{m \times m}$ such that $A V = V Z$. 

For a matrix $A \in \mathbb{S}_n,$ $A >0$ ($ \geq 0, < 0, \leq 0$)  indicates that $A$ is a positive definite (positive semidefinite, negative definite, negative semidefinite) matrix. The sets of $n \times n$ positive (semi-) definite and negative (semi-) definite matrices are denoted by $\mathbb{S}^+_n (\overline{\mathbb{S}}^+_n)$ and $\mathbb{S}^-_n (\overline{\mathbb{S}}^-_n)$, respectively. The notation $\R^-$ ($\R^+$) denotes the set of negative (positive) real numbers, whereas $\C^-$ ($\C^+$) denotes the set of complex numbers with negative (positive) real parts. For a positive integer $q$, the set $\{1,2,\ldots,q\}$ is denoted by $\mathcal{I}_q$. 

We use the shorthand notation $J$ to denote the $2n \times 2n$ skew-symmetric matrix $\begin{bmatrix}
    \mathbf{0}_n & I_n \\ -I_n & \mathbf{0}_n
\end{bmatrix}.$ A matrix $H \in \mathbb{R}^{2n \times 2n}$ is called a \textit{Hamiltonian matrix} if $(J H)^T=(J H).$ A Hamiltonian matrix $H$ has the property that if $\lambda \in \sigma(H)$, then, $-\lambda, \bar{\lambda}, -\bar{\lambda} \in \sigma(H)$. If $(\lambda,v)$ is a right eigenpair of $H$, then, $(-\lambda,(Jv)^T)$ is a left eigenpair of $H$. 

Given matrices $A,Q,D \in \mathbb{R}^{n \times n},$ we let $\mathcal{R}(A,Q,D)$ denote in short the algebraic Riccati equation 
\begin{equation}\label{eqn:ARE}
A^TX + XA + Q - XDX = 0.
\end{equation}
Any matrix $X$ satisfying \eqref{eqn:ARE} is a \textit{solution} to $\mathcal{R}(A,Q,D)$. A solution is called \textit{stabilizing} if the matrix $A-DX$ is Hurwitz. A stabilizing solution, if exists, is unique \cite{zhou1996robust}.

\section{Problem Formulation}\label{sec:Problem Statement}

We consider the problem of privately outsourcing the solution of an algebraic Riccati equation \eqref{eqn:ARE} to a cloud. We make the following assumptions about the cloud:

\begin{assum}\label{a:curious honest adv}
 (i) The cloud is an honest but curious adversary. (ii) Cloud a priori knows the privacy-preserving mechanism employed by the system owner, but does not know the algorithm's \textit{secrets}, i.e. the confidential choices made during the execution of the algorithm.\footnote{See also the first paragraph of Section \ref{sec:Privacy analysis}.}
\end{assum}

The honest but curious assumption means that the cloud follows the protocol instructed by the system owner and
does not deviate from it; however, it stores the information received
from the system to potentially learn more than what has been directly communicated to the cloud. The second item in Assumption \ref{a:curious honest adv} depicts the side knowledge available to the cloud, which asks for the algorithm to be secure even if the cloud knows all the steps of the employed scheme (see also Kerckhoffs’ principle in cryptography \cite{katz2007introduction}). 

The idea that we put forward here is to construct a modified parameter set, namely $(\tilde{A},\tilde{Q},\tilde{D})$, and send it to the cloud such that the solution to $\mathcal{R}(A,Q,D)$ remains the same as the one to $\mathcal{R}(\tilde{A},\tilde{Q},\tilde{D})$ computed by the cloud. 
Such a modification would obviate the need for an algebraic decoding step since the result of the computation can be readily used by the system for control purposes. The focus of this work is on system-theoretic applications, e.g. LQR and $H_{\infty}$ control, where one is interested in the symmetric stabilizing solution of the ARE, i.e., 
$X=X^T$ with $\sigma(A-DX)\subseteq \C^-$. 
To ensure existence, and thus uniqueness \cite{lancaster1995algebraic}, of the stabilizing solution, we impose the following standing assumption:

\begin{assum}\label{a:exist}
    The pair $(A,D)$ is stabilizable and $(Q,A)$ is detectable.
\end{assum}

Before formalizing the problem, we discuss two trivial modifications of the ARE parameters. The first one is to multiply all the coefficients by the same nonzero scalar i.e., $(\tilde{A},\tilde{Q},\tilde{D}) = (kA, kQ, kD)$. However, such scaling barely provides any privacy protection since the system parameters will be revealed up to a scaling factor, e.g. a scaled version of all eigenvalues of the (state) matrix $A$ will be revealed. Moreover, if the cloud suspects/knows the employed scheme, then side knowledge on a single nonzero element of $A$, $Q$, or $D$ is sufficient to identify the secret $k$. 
The second trivial solution is to apply a similarity transformation to change the coordinates, and send the transformed matrices $(T^{-1}AT,T^{T}QT,T^{-1}DT^{-T})$ to the cloud. The obtained solution $\tilde X$ must be then mapped back to a solution of the original ARE by inverting the applied transformation, i.e $X=T^{-T}\tilde{X}T^{-1}$. In this case, however, the eigenvalues of the (system) matrix $A$ remain invariant under similarity transformations and thus will be revealed to the cloud. Furthermore, in system-theoretic applications, a mere change of coordinates does not affect the Markov parameters of the system, which implies that the cloud can completely identify the input-output map, namely the transfer matrix of the associated linear system.
The major privacy breaches mentioned above motivate us to seek for a more sophisticated algebraic perturbation technique for preserving privacy. This gives rise to the following problem formulation:

\begin{prbm}\label{prb:Prob wo realizability}
Let Assumption \ref{a:exist} hold, and let $P$ be the unique stabilizing solution to $\mathcal{R}(A,Q,D).$ Find a nontrivial\footnote{The qualifier ``nontrivial" is meant to exclude the scaled matrices mentioned before. Similarity transformations are automatically excluded from the admissible modifications since the stabilizing solution is different than the one of the original ARE in that case; see also the discussion preceding Problem \ref{prb:Prob wo realizability}.} modified set of parameters $(\tilde{A},\tilde{Q},\tilde{D}) \in \mathbb{R}^{n \times n} \times \mathbb{S}_n \times \mathbb{S}_n$ with $\tilde{A}\neq A,$ $\tilde{Q} \neq Q,$ $\tilde{D} \neq D,$ such that $P$ is the unique stabilizing solution to $\mathcal{R}(\tilde{A},\tilde{Q},\tilde{D}).$
\end{prbm}

While in some applications such as optimal control with conflicting objectives and state feedback for $H_{\infty}$ optimal control \cite{chen1992necessary, bini2011numerical} the matrices $Q$ and $D$ do not need to be sign definite, in other cases both matrices are required to be positive semidefinite; a case in point is given by the LQR problem. 
In the case of LQR, for instance, this implies that the modified Riccati equation parameters $(\tilde{A},\tilde{Q},\tilde{D})$ may not correspond to an LQR problem for any linear system. As such, the cloud realizes that the system parameters have been modified from the original ones. While this can be acceptable, a realizable set of parameters bring additional level of privacy, since the cloud would not even notice that the original parameters have been changed. 
To ensure such realizability, the matrices $\tilde{Q}$ and $\tilde{D}$ must belong to the set of positive semidefinite matrices. This guarantees that $\tilde Q$ is a valid performance matrix for a LQR problem and that $\tilde{D}$ can be written as $\tilde{B} \tilde{R}^{-1} \tilde{B}^T$ for some matrices $\tilde{B}$ and $\tilde{R}>0$. 
This brings us to a stricter version of Problem \ref{prb:Prob wo realizability}:

\begin{prbm}\label{prb:Prob w realizability}
Let $P$ be the unique stabilizing solution to $\mathcal{R}(A,Q,D)$, with $Q,D \in \overline{\mathbb{S}}^+_n$.
Find a nontrivial modified set of parameters $(\tilde{A},\tilde{Q},\tilde{D}) \in \mathbb{R}^{n \times n} \times \overline{\mathbb{S}}^+_n \times \overline{\mathbb{S}}^+_n$ with $\tilde{A}\neq A,$ $\tilde{Q} \neq Q,$ $\tilde{D} \neq D,$ that solves Problem \ref{prb:Prob wo realizability}.
\end{prbm}

As will be observed, a distinctive feature of the proposed approach to solving the two formulated problems is that the perturbed and the original parameter sets are not related via any deterministic maps, which enhances privacy (at the expense of additional computation) and sets it apart from generalized similarity transformations and isomorphisms \cite{sultangazin2020symmetries,xu2015secure,wu2016preserving}.

\section{Proposed eigenvalue shifting technique}\label{sec:Main Results}
Given the triple $(A, Q, D)$, define the map $\mathcal{H}:\mathbb{R}^{n \times n} \times \mathbb{S}_n \times \mathbb{S}_n \to \mathbb{R}^{2n \times 2n}$ as 
\begin{equation*}\label{eqn:Hamiltonian matrix defn}
    \mathcal{H}(A,Q,D) = \begin{bmatrix}
A & -D \\ -Q & -A^T
\end{bmatrix}:=H.
\end{equation*}
The matrix $H$ is the Hamiltonian matrix associated with the algebraic Riccati equation $\mathcal{R}(A,Q,D)$. The Hamiltonian matrix plays a central role in computation of solutions to the ARE \cite{zhou1996robust}. 

We first provide a solution to Problem \ref{prb:Prob wo realizability} and then address Problem \ref{prb:Prob w realizability} in Subsection \ref{subsec:Realizability}.

\subsection{Solution to Problem \ref{prb:Prob wo realizability}}\label{subsec:Solution to Problem 1}
The theorem below provides a solution to Problem \ref{prb:Prob wo realizability}, where the modified set of parameters is obtained by shifting any real eigenvalue of the Hamiltonian matrix. The solution is inspired by the numerical eigenvalue shifting technique in \cite[Sec. 2.4.2]{bini2011numerical}, where the \textit{zero}\footnote{Note that here $H$ does not have any zero or purely imaginary eigenvalues due to
Assumption \ref{a:exist} \cite{lancaster1995algebraic}.} eigenvalue of the Hamiltonian matrix is shifted to improve the convergence rate of the involved computational algorithms.  

\begin{thm}\label{thm:Main theorem mod}
 Let Assumption \ref{a:exist} hold and $P$ be the unique stabilizing solution to $\mathcal{R}(A,Q,D).$ Let $(\lambda_i,v_i) \in \R^- \times \C^{2n}$ be a right eigenpair of $H=\mathcal{H}(A,Q,D)$ where $v_i$ is a unit vector. Then, for any $\Delta \lambda_i \in \R$, there exist matrices $\tilde{A} \in \R^{n \times n}$, $\tilde{Q},\tilde{D} \in \mathbb{S}_{n}$ and correspondingly a Hamiltonian matrix $\tilde{H}:=\mathcal{H}(\tilde{A},\tilde{Q},\tilde{D})$ such that  
    \begin{equation}\label{eqn:Perturbation of H pm 1}
        H +  \Delta \lambda_i \left( v_i p^T_i - q_i q^T_i \right) = \tilde{H}
   \end{equation}
   where \begin{align}\label{eqn:p and q}
     p_i=(J+I_n)v_i,  \quad q_i= Jv_i.
   \end{align}
   Moreover, $P$ is the unique stabilizing solution to $\mathcal{R}(\tilde{A},\tilde{Q},\tilde{D})$ if $\Delta \lambda_i < -\lambda_i$.
\end{thm}
\textit{Proof:} See the appendix for proof.

While Theorem 1 is stated for perturbing real eigenvalues, a more elaborate modification can be used in the case of complex eigenvalues. For brevity, we provide the results only for complex eigenvalues with partial multiplicity $1$. 
\begin{thm}\label{thm:Main theorem cmplx}
    Let Assumption \ref{a:exist} hold and $P$ be the unique stabilizing solution to $\mathcal{R}(A,Q,D).$
    Let $(\mu_j,v_j)$ and $(-\mu_j,v_{-j})$ be two right eigenpairs of $H=\mathcal{H}(A,Q,D)$, where $\mu_j$ ($-\mu_j$) is a complex eigenvalue with partial multiplicity $1$. Then, for any $\Delta \mu_j \in \R,$ there exist matrices $\tilde{A} \in \R^{n \times n}$, $\tilde{Q},\tilde{D} \in \mathbb{S}_{n}$ and correspondingly a Hamiltonian matrix $\tilde{H}:=\mathcal{H}(\tilde{A},\tilde{Q},\tilde{D})$ such that  
   \begin{equation}\label{eqn:Perturbation of H cmplx}
        H +  2 \Delta \mu_j {\rm Re}\left( v_j p^T_j + v_{-j} q^T_j \right) = \tilde{H}
   \end{equation}
   where 
   \begin{align}\label{eqn:p and q cmplx}
        p_j=\theta_j Jv_{-j}, & \qquad q_j=\theta_j Jv_j
   \end{align}
   with the scalar $\theta_j:= (v^T_jJv_{-j})^{-1}$.
   Moreover, $P$ is the unique stabilizing solution to $\mathcal{R}(\tilde{A},\tilde{Q},\tilde{D})$, if $\frac{\Delta \mu_j}{{\rm Re}(\mu_j)}>-1$.
\end{thm}
\textit{Proof:} See the appendix for proof.

The results of Theorem \ref{thm:Main theorem mod} and Theorem \ref{thm:Main theorem cmplx} provide a method to construct a new set of matrices $\tilde{A}, \tilde{Q}$ and $\tilde{D}$ solving Problem \ref{prb:Prob wo realizability}. Notably, the suggested modifications do not spoil the structure of the Hamiltonian matrix, namely \eqref{eqn:Perturbation of H pm 1} and \eqref{eqn:Perturbation of H cmplx} can be viewed again as \textit{Hamiltonian matrices} corresponding to the modified parameters $(\tilde{A},\tilde{Q},\tilde{D})$. 
    
Now, a privacy-preserving outsourcing algorithm can be designed based on Theorems \ref{thm:Main theorem mod} and \ref{thm:Main theorem cmplx}. In particular, we can compute a real eigenpair $(\lambda_i,v_i)$ and then construct a new Hamiltonian matrix using \eqref{eqn:Perturbation of H pm 1}. Similarly, we can compute two complex eigenpairs $(\mu_j,v_j)$ and $(-\mu_j,v_{-j})$ and apply the modifications using \eqref{eqn:Perturbation of H cmplx} in Theorem \ref{thm:Main theorem cmplx}. The new set of matrices obtained from the new Hamiltonian matrix are such that the matrix $P$ is also the unique stabilizing solution to $\mathcal{R}(\tilde{A},\tilde{Q},\tilde{D})$ as desired. Note that $\Delta \lambda_i$ and $\Delta \mu_i$ can be chosen {\em randomly} from the interval dictated by the inequalities stated in Theorems \ref{thm:Main theorem mod} and \ref{thm:Main theorem cmplx}. Finally, we note that, in the same vein, the stated theorems can be iteratively applied to the new Hamiltonian matrix, thus modifying additional eigenvalues before submitting the query (modified ARE) to the cloud.     

\begin{rem}
Numerical methods available for solving algebraic Riccati equations have a time complexity of at least $\mathcal{O}(n^3)$ for a system of size $n$ \cite{bini2011numerical}. In comparison, a privacy preserving outsourcing algorithm that modifies $2k$ eigenvalues of Hamiltonian matrix based on Theorems \ref{thm:Main theorem mod} and \ref{thm:Main theorem cmplx} has a time complexity of $\mathcal{O}(k n^2)$, which is a large difference for large-scale systems with $n>>k$. As $k$ increases, the computational load of the algorithm increases proportionally. At the same time, the discrepancy between the modified parameters and the original ones increases as more eigenvalues are perturbed, thereby enhancing privacy. This leads to a trade-off between local computational load and privacy guarantees.
\end{rem}

\subsection{Solution to Problem \ref{prb:Prob w realizability} }\label{subsec:Realizability}
In this subsection, we shift our attention to Problem 2, which asks for preserving positive semidefiniteness of $Q$ and $D$ in the employed perturbation. As previously mentioned, this would prevent the cloud from knowing that the parameters have been perturbed. This problem is more complex than the former one, and to tame the complexity we restrict our perturbation method to real eigenvalues. 

We decompose the vectors $v_i$, $p_i$ and $q_i$ defined in Theorem \ref{thm:Main theorem mod} as $v_i={\rm col}(v_{iu},v_{il})$, $p_i={\rm col}(p_{iu},p_{il})$, and $q_i={\rm col}(q_{iu},q_{il})$,
where $v_{iu}, v_{il},p_{iu},p_{il},q_{iu},q_{il} \in \R^n$.
As a result of perturbing $H$ to $\tilde{H}$ using \eqref{eqn:Perturbation of H pm 1}, we have
\begin{align}\label{eqn:Changes in D and Q}
    \tilde{D} = D + (\Delta\lambda_i)F_i, & \qquad
    \tilde{Q} = Q + (\Delta\lambda_i)G_i, 
\end{align}
where $F_i:=q_{iu} q^T_{il} - v_{iu} p^T_{il}$ and $G_i:=q_{il} q^T_{iu} - v_{il} p^T_{iu}$. Note that since $\tilde{H}$ is a Hamiltonian matrix by Theorem \ref{thm:Main theorem mod}, the matrices $F_i$ and $G_i$ are symmetric.
If $F_i$ and $G_i$ are either both positive semidefinite or negative semidefinite, then we can pick the sign of the perturbation $\Delta \lambda_i$ accordingly such that $\tilde{D}$ and $\tilde{Q}$ remain positive semidefinite and thus solve Problem \ref{prb:Prob w realizability}. 
As we rely on local computations to check the definiteness of such matrices, we state the following algebraic result which reduces 
this check to verifying the sign of scalar quantities. For convenience, we drop the index $i$, as the forthcoming expressions remain the same for different eigenvalues.
\begin{prop}\label{prop:Semidefiniteness of L}
Define the scalar quantities 
\begin{equation}\label{eqn:algebraic quantities for realizability}
    \begin{aligned}
        \alpha^F_{uu} := q^T_{u}q_{u}, &~~~~ \beta^F_{uu} := v^T_{u}q_u, &~~~~ \gamma^F := q^T_up_l,
         \\ \alpha^F_{ul} := q^T_uq_l, &~~~~ \beta^F_{ul} := v^T_uq_l, &~~~~
        \delta^F := v^T_up_l .
    \end{aligned}
\end{equation}
Then, it holds that
\begin{align*}
        F \geq 0, & \text{ if} & \alpha^F_{ul} \delta^F - \gamma^F \beta^F_{ul}  \leq 0 & \text{ and} & \alpha^F_{uu} \alpha^F_{ul} - \beta^F_{uu} \gamma^F \geq 0 \\
        F \leq 0, & \text{ if} & \alpha^F_{ul} \delta^F - \gamma^F \beta^F_{ul}  \leq 0 & \text{ and} & \alpha^F_{uu} \alpha^F_{ul} - \beta^F_{uu} \gamma^F \leq 0.
\end{align*}
\end{prop}
\textit{Proof:} See the appendix for proof.

 If $F \geq 0\, (\leq 0)$, then choosing $\Delta \lambda_i \geq 0\, (\leq 0)$ ensures that $\tilde{D} \geq 0$; see \eqref{eqn:Changes in D and Q}.
Analogously, Proposition \ref{prop:Semidefiniteness of L} can be stated to verify positive/negative definiteness of $G$ by replacing the superscripts $F$ with $G$ in the inequalities and defining 
\begin{equation}\label{eqn:algebraic quantities for realizability 2}
    \begin{aligned}
        \alpha^G_{ll} := q^T_{l}q_{l}, &~~~~ \beta^G_{ll} := v^T_{l}q_l, &~~~~ \gamma^G := q^T_lp_u,
         \\ \alpha^G_{lu} := q^T_lq_u, &~~~~ \beta^G_{lu} := v^T_lq_u, &~~~~
        \delta^G := v^T_lp_u.
    \end{aligned}
\end{equation}

Other than the conditions mentioned in Proposition \ref{prop:Semidefiniteness of L}, we mention another scenario where positive semidefiniteness of $\tilde{Q}$ and $\tilde{D}$ can be guaranteed. 
This concerns the case where ${\rm ker}(D) \subseteq {\rm ker}(F)$ and ${\rm ker}(Q) \subseteq {\rm ker}(G)$.
Denoting the minimum (maximum) nonzero eigenvalue of a matrix $A$ by $\mu_{min}(A)$ ($\mu_{max}(A)$), the following result addresses the perturbation guideline for this case, which is particularly useful in the case when $F$ is sign indefinite. 

\begin{prop}\label{prop:Image of L and D}
Assume that ${\rm ker}(D) \subseteq {\rm ker}(F)$.
Then, 
    \begin{equation}\label{eqn:Lambda range for positive tilde D}
        \Delta \lambda \in [\underline{\lambda}^F,\overline{\lambda}^F] \implies \tilde{D} \geq 0,
    \end{equation}
    where 
    \begin{equation}\label{eqn:Upper and lower bounds for L}
        \underline{\lambda}^F = -\frac{\mu_{min}(D)}{\mu_{max}(F)}, \quad
        \overline{\lambda}^F = -\frac{\mu_{min}(D)}{\mu_{min}(F)}.
    \end{equation}
\end{prop}
\textit{Proof:} See the appendix for proof.

Proposition \ref{prop:Image of L and D} analogously applies to $\tilde{Q}$ as well, by replacing $D$ with $Q$ and the superscript $F$ with $G$. 
Subsequently, a suitable perturbation ensuring $\tilde{D} \geq 0$ and $\tilde{Q} \geq 0,$ can be chosen as 
$\Delta \lambda_i \in [\lambda^l_i,\lambda^u_i]$, with
\begin{align}\label{eqn:Upper and lower bounds}
    \lambda^l = \max (\underline{\lambda}^F,\underline{\lambda}^G),& \qquad \lambda^u = \min (\overline{\lambda}^F,\overline{\lambda}^G).
\end{align}
We note that, unlike Proposition \ref{prop:Semidefiniteness of L},  Proposition \ref{prop:Image of L and D} can be applied to the cases where the matrix $F$ (analogously $G$) is not sign definite. On the other hand, it asks for additional assumption that ${\rm ker}(D) \subseteq {\rm ker}(F)$ (analogously, ${\rm ker}(Q) \subseteq {\rm ker}(G)$). Noting that ${\rm{rank}}(F)\leq 2$ (${\rm{rank}}(G)\leq 2$), the latter condition is likely to be satisfied in the case where the rank of matrix $D$ ($Q$) is very high. 

The two scenarios discussed in Propositions \ref{prop:Semidefiniteness of L} and \ref{prop:Image of L and D} can be summarized into Algorithm \ref{alg:Algorithm 2} below which solves Problem \ref{prb:Prob w realizability} under the stated hypotheses. 

\begin{savenotes}
\begin{algorithm}
{\footnotesize{
\caption{Algorithm to solve Problem \ref{prb:Prob w realizability}}\label{alg:Algorithm 2}
\textbf{Input:} Hamiltonian matrix $H$\\
    \textbf{for $i=1,2,\ldots, n$}
    \begin{enumerate}
        \item Compute a real eigenvalue $\lambda_i < 0$.
        \item Compute the corresponding right eigenvector $v_{i},$ and construct $p_i$ and $q_i$ from it using \eqref{eqn:p and q}. 
       \item
        \textbf{If } $F_i \geq 0$ and $G_i \geq 0$\\ Choose a $\Delta \lambda_i \in (0,-\lambda_i)$, and \textbf{go to} Step 5. \\
            \textbf{else if }\footnote{A sufficient condition for positive/negative semidefiniteness of $F_i$ (similarly, $G_i$) is provided in Proposition \ref{prop:Semidefiniteness of L}.} $F_i \leq 0$ and $G_i \leq 0$\\ Choose a $\Delta \lambda_i \in (-\infty,0)$,and \textbf{go to} Step 5. \\
            \textbf{else }continue \\
            \textbf{end}
        \item \textbf{If ${\rm ker}(D) \subseteq {\rm ker}(F_i)$ and ${\rm ker}(Q) \subseteq {\rm ker}(G_i)$}
        \begin{itemize}
            \item Compute $\lambda^u_i$ and $\lambda^l_i$ using \eqref{eqn:Upper and lower bounds for L} and \eqref{eqn:Upper and lower bounds}.
            \item Choose a $ \Delta \lambda_i \in (\lambda^l_i, \min(-\lambda_i, \lambda^u_i))$, and \textbf{go to} Step 5.
        \end{itemize}
        \textbf{else skip} to the next iteration.\\
        \textbf{end}
        \item Construct $\tilde{H}$ using \eqref{eqn:Perturbation of H pm 1}.\\ 
        \textbf{stop}
    \end{enumerate}
    \textbf{end}\\
\textbf{Output:} Modified Hamiltonian matrix $\tilde{H}.$
}}
\end{algorithm}
\end{savenotes}

\begin{rem}\label{rem:Remark on ensuring stabilizing solution}
    The bounds on $\Delta \lambda_i$ in Step 3 (resp. Step 4) of Algorithm \ref{alg:Algorithm 2} are determined by intersecting the conditions in Theorem 1 (i.e., $\Delta \lambda_i < -\lambda_i$) with those resulting from Proposition \ref{prop:Semidefiniteness of L} (resp. Proposition \ref{prop:Image of L and D}). 
\end{rem}

\section{Privacy analysis}\label{sec:Privacy analysis}
We perform the following privacy analysis with a specific focus on the Hamiltonian matrix.\footnote{For brevity, we only provide the expressions for the privacy preserving mechanism of Problem \ref{prb:Prob wo realizability}, with {real} eigenvalue perturbations.} The significance of the Hamiltonian matrix lies in the fact that it captures the combined information of the system dynamics and optimization parameters. Moreover, in view of Assumption \ref{a:curious honest adv}, cloud is aware that the privacy preserving mechanism is based on the perturbation of the Hamiltonian, and thus, as a curious party, seeks to reconstruct the true Hamiltonian matrix, denoted in this section by $H^\star$, from the received modified one $\tilde{H}$.
The uncertainty that the cloud is facing in such a reconstruction originates from the algorithms \textit{secrets}, namely the choices 
(henceforth called the \textit{indices}) of the perturbed eigenvalues and the amount of perturbations (henceforth
called the \textit{magnitudes}, with a slight abuse of the terminology).
In case of perturbing a single real eigenvalue, the confusion set faced by the cloud can be explicitly given by 
\begin{multline}\label{eqn:Confusion set k=1}
    \mathcal{C}_1 =\{H|H=\tilde{H}-\gamma f(i,\tilde{H}), i \in \mathcal{I}_{r}, \gamma > \lambda_{i}(\tilde{H})\},
\end{multline}
where $r$ is the number of negative real eigenvalues of $\tilde{H}$, $\mathcal{I}_r$ is defined in Section \ref{sec:Preliminaries}, and $\lambda_{i}(\tilde{H})$ denotes the $i$th negative real eigenvalue of $\tilde{H}$. The mapping $f(i,\tilde{H})$ outputs the matrix $v_i p^T_i - q_i q^T_i$, where $(\lambda_i,v_i)$ denotes the $i$th eigenpair of the perturbed Hamiltonian matrix $\tilde{H}$ and the vectors $p_i$ and $q_i$ are constructed from $v_i$ according to \eqref{eqn:p and q}. The cloud's knowledge that $\lambda_i(H^\star)<0$ establishes the lower bound on the magnitude $\gamma$ in \eqref{eqn:Confusion set k=1}. Note that any $H\in \mathcal{C}_1$ is consistent with the knowledge of the cloud on the privacy mechanism and the received perturbed Hamiltonian $\tilde{H}$. To reconstruct the true Hamiltonian, the cloud needs to correctly guess the index $i$ from the set $\mathcal{I}_{r}$ of cardinality $r$. If the index is correctly guessed, then the confusion set will be parameterized by the single scalar parameter $\gamma$. Therefore, the extent of ambiguity in finding $i$ and $\gamma$ can be succinctly captured by the pair $(r, 1)$, serving as a privacy measure in this special case.

For the case where $k>1$ eigenvalues of the Hamiltonian are perturbed, to reconstruct the true Hamiltonian, the cloud needs to guess a series of $k$ indices and $k$ magnitudes. To generalize \eqref{eqn:Confusion set k=1}, i.e. writing the confusion set faced by the cloud, we recursively define the following sets for each $j=0,1,\ldots, k-1$:
\begin{multline*}
    \mathcal{C}_{j+1} := \{H|H=H_j-\gamma_j f(i_j,H_j), H_j \in \mathcal{C}_j,\\ i_j \in \mathcal{I}_q, \gamma_j > \lambda_{i_j}(H_j)\},
\end{multline*}
where $q=r-j$ and $\mathcal{C}_0:=\{\tilde{H}\}.$ 
The final set, namely $\mathcal{C}_k$ depicts the confusion set faced by the cloud, which generalizes \eqref{eqn:Confusion set k=1} to the case $k>1$. 
Bearing in mind that $\gamma_j=0$ is an admissible value in the above set, we have that $\mathcal{C}_j \subseteq \mathcal{C}_{j+1}$ for each $j$.
Therefore, we obtain 
 $   \{\tilde{H}\}=\mathcal{C}_0 \subseteq \mathcal{C}_1 \subseteq \cdots \subseteq \mathcal{C}_{k-1} \subseteq \mathcal{C}_k,
$
establishing the observation that the confusion set keeps growing with every additional perturbation. To see how the privacy measure can be extended to the case $k>1$, we note that cloud has to correctly guess the correct $k-$long sequence of indices $(i_0, \ldots, i_{k-1})$ out of $\frac{r!}{(r-k)!}$ possible choices. If the sequence is correctly guessed, the confusion set will be parameterized by $k$ scalar parameters $\gamma_j$, with $j=0, 1, \ldots, k-1$. Consequently, the extent of the ambiguity faced by the cloud can be compactly represented by the pair $\big(\frac{r!}{(r-k)!},k\big)$ in this case.

\section{Case study}\label{sec:Case Study}
We illustrate the algorithm by applying it to the LQR problem. We work with the benchmark for ARE in \cite{abels1999carex}.  
This corresponds to the LQR problem \cite[Ch. 6]{liberzon2011calculus} for a linear dynamical system
$
\dot{x}(t)=Ax(t)+Bu(t), \; y(t)=Cx(t),
$
of size $n=100$ with $m=10$ control inputs and $p=10$ outputs. The system matrices can be chosen such that
$Q=CC^T$ and $D=BR^{-1}B^T$, where the performance matrix $R$ can be taken as the identity matrix for simplicity. 
We run the privacy-preserving outsourcing algorithm to modify $2r=20$ real eigenvalues. 

The norms of the matrices $\tilde{A}-A$, $\tilde{D}-D$, and $\tilde{Q}-Q$ quantify how different the modified parameters are from the original ones. Therefore, the expressions $\frac{||\tilde{A}-A||}{||A||}$, $\frac{||\tilde{D}-D||}{||D||}$, and $\frac{||\tilde{Q}-Q||}{||Q||}$ can serve as measures of privacy for the respective matrices. Table \ref{tbl:Table of privacy metrics} shows the behaviour of these privacy measures with the number of iterations. 
\begin{table}[h!]
\vspace{2mm}
\centering
    \begin{tabular}{|c||c|c|c|}
         \hline
         Iteration & $\frac{||\tilde{A}-A||}{||A||}$ & $\frac{||\tilde{D}-D||}{||D||}$ & $\frac{||\tilde{Q}-Q||}{||Q||}$ \\
         \hline\hline
         1 & 0.1917 & 0.1955 & 0.2082 \\
         \hline
         5 & 0.3024 & 0.3291 & 0.3148 \\
         \hline
         9 & 0.3957 & 0.3841 & 0.3883 \\
         \hline
    \end{tabular}    
    \caption{Evolution of privacy measures with number of iterations of the privacy preserving algorithm}
\label{tbl:Table of privacy metrics}
\end{table}
For comparison, we run Algorithm \ref{alg:Algorithm 2} for the same system only for {\emph{one iteration}}, i.e. perturbing only one eigenvalue. The first row of Table \ref{tbl:Table of privacy metrics} in this case reduces to $0.0151$, $0.0127$, and $0.0157$ respectively, for the matrices $A$, $D$, and $Q$. 
These values are much smaller 
since a narrower range of perturbations is allowed in Algorithm \ref{alg:Algorithm 2} due to the additional realizability requirement. 
\section{Conclusions}\label{sec:Conclusions}
We have addressed the issue of securely outsourcing the solution of algebraic Riccati equation (ARE) to a cloud. We formulated the corresponding privacy-preserving problem with and without a realizability constraint, which determines whether the perturbed parameters must preserve their semidefiniteness. The proposed algorithms provided recipes for perturbing the eigenvalues of the Hamiltonian matrix in such a way that the solution to the corresponding ARE remains unchanged and thus usable to compute the actual solution. Finally, we have illustrated our approach through a numerical example for a linear quadratic regulator (LQR). 
It is of interest to explore other control problems where the same principles can be applied, namely where the demanding part of the computation can be performed on dummy system models in a way that the computed result can still be used to derive the desired controller for the true system.

\section{Appendix: Proofs of the results}\label{sec:Appendix}
The following result can be readily obtained from Rado's theorem \cite{bru2012brauer} and 
\cite[Theorem 1.11]{bini2011numerical}:
\begin{lem}\label{lem:Rado's thm}
    Let $L \in \C^{n \times n}$ and $\sigma(L)=\{\lambda_1,\ldots,\lambda_n\}$. Let $M \in \C^{n \times r}$ be a matrix of full column rank whose columns are right eigenvectors of $L$ corresponding to the eigenvalues $\lambda_1,\ldots,\lambda_r$. Denote $\Lambda={\rm diag}(\lambda_1,\ldots,\lambda_r)$. For any matrix $N \in \C^{r \times n}$, the following statements hold:
    \begin{enumerate}
        \item $\sigma(L+MN)=\sigma(\Lambda+NM) \cup \{\lambda_{r+1},\ldots,\lambda_n\}$.
        \item Let $V$ be a matrix of full column rank such that $LV=VZ$ and $M=V\tilde{M}$ for some matrices $Z$ and $\tilde{M}$. Then, $(L+MN)V=V(Z+\tilde{M}\tilde{N})$ where $\tilde{N}=NV$. In addition, $Z \tilde{M}= \tilde{M}\Lambda$ and $\sigma(Z+\tilde{M}\tilde{N})=\left(\sigma(Z)\backslash\sigma(\Lambda)\right) \cup \sigma(\Lambda+NM)$. 
        \item Let $W$ be a matrix whose rows are left eigenvectors of $L$ corresponding to eigenvalues $\sigma(L)\backslash \sigma(\Lambda)$. Then, $W(L+MN) = WL$.
    \end{enumerate}
\end{lem}
We note that analogous results can be stated for perturbations using the left eigenvectors of $L$.  

\textit{Proof of Theorem \ref{thm:Main theorem cmplx}}:
\footnote{For convenience, we postpone the proof of Theorem \ref{thm:Main theorem mod}.}
Since matrix $P$ is the stabilizing solution to $\mathcal{R}(A,Q,D)$, we have that  \cite{lancaster1995algebraic,zhou1996robust} 
\begin{equation}\label{eqn:H-invariance}
    \begin{aligned}
        H {\rm col}(I_n,P) &= {\rm col}(I_n,P) (A-DP),\\
        \sigma(A-DP)&=\sigma(H) \cap \C^-
    \end{aligned}
\end{equation}
In order to prove Theorem \ref{thm:Main theorem cmplx}, it is sufficient to show that 
\begin{itemize}
    \item the perturbation in \eqref{eqn:Perturbation of H cmplx} is a Hamiltonian matrix,
    \item the subspace ${\rm im}\left({\rm col}(I_n,P)\right)$ is $\tilde{H}-$invariant corresponding to the eigenvalues $\sigma(\tilde{H}) \cap \C^-$.
\end{itemize}

 To this end, we represent the perturbation in \eqref{eqn:Perturbation of H cmplx} as the matrix product $MN$ where $M:= \begin{bmatrix}
     v_{j} & v_{-j} & \bar{v}_j & \bar{v}_{-j}
 \end{bmatrix}$ and $N:=\Delta \mu_j \begin{bmatrix}
     p_j & q_j & \bar{p}_j & \bar{q}_j 
 \end{bmatrix}^T$. 
  It can be easily seen that
 \begin{equation*}
     MN=\Delta \mu_j\theta_j\left( v_jv^T_{-j} + v_{-j}v^T_{-j} + \bar{v}_j\bar{v}^T_{-j} + \bar{v}_{-j}\bar{v}^T_j \right)J^T
 \end{equation*}
 is a Hamiltonian matrix, as it satisfies the equation $(JMN)^T=JMN$.
 Noting that $(\mu_j,p^T_j)$, $(-\mu_j,q^T_j)$, $(\bar{\mu}_j,\bar{p}^T_j)$ and $(-\bar{\mu}_j,\bar{q}^T_j)$ are left eigenpairs of $H$, and that the partial multiplicity of $\mu_j$ is assumed to be $1$, we have that the matrix $NM$ is diagonal:
 \begin{align*}
     NM&=\Delta \mu_j{\rm diag} (p^T_jv_j, q^T_jv_{-j},\bar{p}^T_j\bar{v}_j,\bar{q}^T_j\bar{v}_{-j})\\
     &={\rm diag}(\Delta \mu_j,-\Delta \mu_j,\Delta \mu_j,-\Delta \mu_j).
 \end{align*}
 By using Lemma \ref{lem:Rado's thm}, we obtain the eigenvalues of the new Hamiltonian matrix $\tilde{H}=H+MN$ as $\sigma(\tilde{H})=(\sigma(H)\backslash\sigma(\Lambda)) \cup \sigma(\Lambda+NM)$. That is, the eigenvalues $\mu_j,-\mu_j,\bar{\mu}_j,-\bar{\mu}_j$ of $H$ are respectively replaced by  $\mu_j+\Delta \mu_j,-\mu_j-\Delta \mu_j,\bar{\mu}_j+\Delta \mu_j,-\bar{\mu}_j-\Delta \mu_j$.
 By using equation \eqref{eqn:H-invariance} and choosing the matrices $V$ and $Z$ of Lemma \ref{lem:Rado's thm} as $V:={\rm col}(I_n,P)$ and $Z:=A-DP$, we find that ${\rm im}(V)$ is also $\tilde{H}-$invariant, corresponding to the eigenvalues $\sigma(Z+\tilde{M}\tilde{N})=\sigma(\tilde{A}-\tilde{D}P)$. As the condition $\frac{\Delta \mu_i}{{\rm Re}(\mu_i)}>-1$ ensures that the signs of real parts of the modified eigenvalues remain unchanged, the subspace ${\rm im}(V)={\rm im}\left({\rm col}(I_n,P)\right)$ is $\tilde{H}-$invariant with $\sigma(\tilde{A}-\tilde{D}P)=\sigma(\tilde{H}) \cap \C^-$, which completes the proof.

\textit{Proof of Theorem \ref{thm:Main theorem mod}}:
Let $M_1:=v_i$, $N_1:=\Delta \lambda_i p^T_i$, $M_2:=q^T_i$, and $N_2:=-\Delta \lambda_i q_i$ with $p_i$ and $q_i$ as defined in \eqref{eqn:p and q}, and rewrite \eqref{eqn:Perturbation of H pm 1} as $\tilde{H}=H+M_1N_1+N_2M_2$. From the assumption that $v_i$ is a unit vector, it follows that $N_1M_1=\Delta \lambda_i p^T_iv_i=\Delta \lambda_i$, and $M_2N_2=-\Delta \lambda_i q^T_iq_i=-\Delta \lambda_i$. Using Lemma \ref{lem:Rado's thm} with $L=H$, we have $\sigma(H+M_1N_1)=\left(\sigma(H)\backslash\{\lambda_i\}\right) \cup \{\lambda_i+\Delta \lambda_i\}$. Notice that due to statement $3$ of Lemma \ref{lem:Rado's thm}, $(-\lambda_i,q^T_i)$ is a left eigenpair of $H+M_1N_1$. Then, by using a modified version of Lemma \ref{lem:Rado's thm} for perturbations using left eigenvectors of $L=H+N_1M_1$, we have that 
\begin{align*}
    \sigma(\tilde{H})=\left(\sigma(H)\backslash \{\lambda_i,-\lambda_i\} \right) \cup \{\lambda_i+\Delta \lambda_i, -\lambda_i-\Delta \lambda_i\}.
\end{align*}
The matrix $\tilde{H}=H+ M_1N_1+N_2M_2$ is a Hamiltonian matrix as $(JM_1N_1+JN_2M_2)^T=JM_1N_1+JN_2M_2$. The fact that $P$ is a stabilizing solution of the modified ARE follows from arguments similar to those in the proof of Theorem \ref{thm:Main theorem cmplx}.

\textit{Proof of Proposition \ref{prop:Semidefiniteness of L}}:
Since $F$ is a symmetric matrix, for any $x \in \R^n$, there exist $x_1 \in {\rm ker}(F)$ and $x_2 \in {\rm im}(F)$ such that $x=x_1+x_2$ and $x^T_1 x_2=0$. Thereby, we obtain $x^TFx=x^T_2 F x_2$. Noting that ${\rm im}(F) = {\rm span}(q_u,v_u)$, we have that if $x^T F x \geq 0 (\leq 0) ~ \forall x \in {\rm span}(q_u,v_u)$, then, $F \geq 0 (\leq 0)$.
Now, let $x = k_1 q_u + k_2 v_u$ for some $k_1,k_2 \in \mathbb{R}.$ Then, $x^T F x = \mathcal{A} k^2_1 + \mathcal{B} k_1k_2 + \mathcal{C} k^2_2$,
where $\mathcal{A}:= \alpha_{uu} \alpha_{ul} - \beta_{uu} \gamma,$ $\mathcal{B}:= \alpha_{uu} \beta_{ul} - \beta_{uu} \delta + \beta_{uu} \alpha_{ul} - \xi \gamma,$ $\mathcal{C}:= \beta_{uu} \beta_{ul} - \xi \delta$ (the superscript $F$ is dropped), and $\xi:=v^T_u v_u$.
If $k_2 = 0,$ then $x^TFx=\mathcal{A}k^2_1$ and its sign is the same as that of $\mathcal{A}.$
If $k_2 \neq 0,$ then $x^TF x = k^2_2 \left(\mathcal{A} (\frac{k_1}{k_2})^2 + \mathcal{B} (\frac{k_1}{k_2}) + \mathcal{C}\right).$ The sign of $x^T F x$ remains the same for all the values of $k_1$ and $k_2$ if and only if the discriminant $\Pi:= \mathcal{B}^2 - 4 \mathcal{A} \mathcal{C} \leq 0$.
After some simplifications, we arrive at $\Pi = 4 (\alpha_{uu} \xi- \beta^2_1) (\alpha_{ul} \delta - \gamma \beta_{ul})$.
From the Cauchy-Schwarz inequality, we find that $\alpha_{uu} \xi - \beta^2_1 = ((q^T_uq_u) (v^T_u v_u) - (v^T_u q_u)^2) \geq 0$. 
Now if $\alpha_{ul} \delta - \gamma \beta_{ul} \leq 0,$ we have $\Pi \leq 0$, which implies that the sign of $x^TFx$ is the same as that of $\mathcal{A}$. This completes the proof.

\textit{Proof of Proposition \ref{prop:Image of L and D}}:
For any $x \in \R^n$, we have 
\begin{equation}\label{eqn:Quadratic form of tilde D}
    x^T\tilde{D}x = x^TDx + \Delta \lambda x^T F x
\end{equation}
from \eqref{eqn:Changes in D and Q}. If $x \in {\rm ker}(F)$, we have $x^T\tilde{D}x = x^TDx \geq 0$.
Now assume that $x \notin {\rm ker}(F)$. Then, $x \notin {\rm ker}(D)$ by the assumption. Since $D$ is symmetric, there exist $x_1 \in {\rm ker}(D) \subseteq {\rm ker}(F)$ and $x_2 \in {\rm im}(D)$ such that $x=x_1+x_2$ and $x^T_1x_2=0$. It follows that $x^TDx=x^T_2Dx_2$ and $x^TFx=x^T_2Fx_2$, which simplifies \eqref{eqn:Quadratic form of tilde D} to
\begin{equation}\label{eqn:Modified QF of tilde D}
    x^T\tilde{D}x = x^T_2Dx_2 + \Delta \lambda x^T_2 F x_2.
\end{equation}
The expression $x^T_2Fx_2$ is bounded as $\lambda_{min}(F)||x_2||^2 \leq x^T_2 F x_2 \leq \lambda_{max}(F)||x_2||^2$, where $\lambda_{min}(\lambda_{max})$ denotes the minimum (maximum) eigenvalue. Moreover, since $x_2 \in {\rm im}(D)$, it holds that $\mu_{min}(D)||x_2||^2 \leq x^T_2 D x_2$,
where $||.||$ stands for the Euclidean norm. Substituting the above inequalities in \eqref{eqn:Modified QF of tilde D}, we arrive at 
\begin{equation*}\label{eqn:Final bounds on QF of tilde D}
    x^T\tilde{D}x \geq  \begin{cases}
        \mu_{min}(D)||x_2||^2 + (\Delta \lambda) \lambda_{min}(F)||x_2||^2, \Delta \lambda > 0, \\
        \mu_{min}(D)||x_2||^2 + (\Delta \lambda) \lambda_{max}(F)||x_2||^2, \Delta \lambda < 0.
    \end{cases}
\end{equation*}
We obtain \eqref{eqn:Lambda range for positive tilde D} and \eqref{eqn:Upper and lower bounds for L} by analyzing the range of values for which the right hand side of the above inequality is nonnegative.

\bibliographystyle{ieeetr}
\bibliography{References}

\end{document}